\documentclass[reqno]{amsart}
\usepackage{verbatim}
\usepackage{hyperref}
\usepackage{comment}
\usepackage{enumerate}
\usepackage{graphicx}
\usepackage{hhline}
\usepackage{amssymb}
\usepackage{graphics,verbatim,color,amsthm}
\usepackage[all]{xy}
\usepackage{fullpage}
\usepackage{tikz}
\usepackage{chngcntr}
\counterwithout{figure}{section}
\counterwithout{figure}{subsection}
\counterwithout{figure}{subsubsection}

\theoremstyle{plain}
\newtheorem{theorem}{Theorem}

\newtheorem{definition}{Definition}

\theoremstyle{remark}

\theoremstyle{definition}
\newtheorem{example}{Example}

\def\sspan{\operatorname{span}}

\begin{document}

\title[]{Complete spelling rules for the Monster tower over three-space}

\author{Alex Castro, Wyatt Howard, and Corey Shanbrom}

\address{Mathematics Department, PUC-Rio,
Rua Marques de Sao Vicente, 225, Edificio Cardeal Leme, Sala 862, Gavea, Rio de Janeiro, CEP 22451-900, Brazil}

\email{alex.castro@mat.puc-rio.br}

\address{
Physical Sciences, Mathematics, and Engineering Division,
De Anza College,
21250 Stevens Creek Blvd,
Cupertino, CA 95014, USA}

\email{howardwyatt@fhda.edu}

\address{Department of Mathematics and Statistics \\ California State University, Sacramento \\ Sacramento, CA 95819, USA}

\email{corey.shanbrom@csus.edu}

\date{\today}



\begin{abstract}  
The Monster tower, also known as the Semple tower, is a sequence of manifolds with distributions of interest to both differential and algebraic geometers.  Each manifold is a projective bundle over the previous.  Moreover, each level is a fiber compactified jet bundle equipped with an action of finite jets of the diffeomorphism group.  There is a correspondence between points in the tower and curves in the base manifold.  These points admit a stratification which can be encoded by a word called the RVT code.  Here, we derive the spelling rules for these words in the case of a three dimensional base.  That is, we determine precisely which words are realized by points in the tower.  To this end, we study the incidence relations between certain subtowers, called Baby Monsters, and present a general method for determining the level at which each Baby Monster is born.  Here, we focus on the case where the base manifold is three dimensional, but all the methods presented generalize to bases of arbitrary dimension.
 \end{abstract}

\maketitle


\section{Introduction}

\subsection{Motivation} 
The Monster tower, also known as the Semple tower, lies in the intersection of differential geometry, non-holonomic mechanics, singularity theory, and algebraic geometry.  Cartan (\cite{C} studied the diffeomorphism group action on jet spaces, which led to developments in the fields of Goursat distributions and sub-Riemannian geometry.  Jean \cite{J}, Luca and Risler \cite{LR}, Li and Respondek \cite{respondek}, Pelletier and Slayman \cite{pelletier1, pelletier2}, and others have studied models of various kinematic systems (a car pulling $n$ trailers, motion of an articulated arm, $n$-bar systems).  Montgomery and Zhitomirskii \cite{MZ1} studied the relationship with curve singularities; later, so did we \cite{castro3, Shanbrom}. And we discovered in \cite{castro2} that algebraic geometers have long studied these objects under different names.  We have begun pursuing these connections \cite{bridges} and working with algebraic geometers to consolidate understanding and improve existing terminology and techniques \cite{kennedy}.  Here, we study the RVT code for the tower, which is invariant under the action of the diffeomorphism group.  This is related to work on the classification problem studied by Mormul \cite{Mo1, Mo3, Mo2}, Montgomery and Zhitomirskii \cite{MZ2, MZ1}, the authors \cite{castro3}, and others.

In the geometric theory of differential equations, we speculate that there may be some interesting connections between the singularity theory of the Monster tower and the general Monge problem for underdetermined systems of ordinary differential equations with an arbitrary number of degrees of freedom.
In \cite{kumpera}, the authors derive sufficient conditions, in terms of truncated multi-flag systems, for the existence of a Monge-Cartan parametrization of the general solution of such systems in the regular case. To our knowledge, no connection has been made with the singular theory of multi-flags presented in this note.  Similar undetermined systems of ordinary differential equations are common in geometric control theory when studying flat outputs of nonlinear control systems \cite{respondek2}.
A detailed account of the geometry of differential equations in jet spaces can be found in \cite{kushner}, where symmetry methods from contact and symplectic geometry are used to solve non-trivial nonlinear partial and ordinary differential equations. 

It remains to investigate the correspondence between finite jets of spatial curves and normal forms of special multi-flags. One should explore the depth of the correspondence between Arnold's A-D-E classification \cite{arnold} and the listing of normal forms of Goursat multi-flags.

Finally, current work with algebraic geometers \cite{kennedy} extends and generalizes the results of this paper to the case of an $n$-dimensional base.  An interesting open question here concerns the existence of moduli in orbits of the action of the diffeomorphism group of the base space.

Thus, it is apparent that this object is of interest to a variety of pure and applied mathematicians, and that it presents a wealth of interesting problems which have potential to shed light in surprising areas.

\subsection{History}
The subject begins with the study of Goursat distributions, which are bracket-generating (completely non-holonomic) but slow growing.
Cartan \cite{C} studied the model of the canonical contact distribution on the jet space $J^{k}(\mathbb{R}, \mathbb{R})$. All Goursat distributions were believed to be equivalent to Cartan's under the action of the diffeomorphism group until Giaro, Kumpera, and Ruiz discovered the first singularity in 1978 \cite{GKR}. 

Jean \cite{J} studied the kinematic model of a car pulling $N$ trailers, a system which is locally universal for Goursat distributions of corank $N+1$.  He developed a geometric stratification given by regions in the configuration space of the model in terms of critical angles.  Montgomery and Zhitomirskii \cite{MZ2} introduced the Monster tower, a sequence of manifolds with distributions in which every Goursat germ occurs, allowing for Jean's strata to be recast in terms of positions of members of a canonical subflag of the Goursat flag.  Mormul \cite{Mo1} labelled the strata from \cite{MZ2} by words in the letters GST, which became the RVT code in \cite{MZ1}.  
  In \cite{MZ1}, Montgomery and Zhitomirskii showed that Goursat germs correspond to finite jets of Legendrian curve germs, and that the RVT coding corresponds to several classical invariants in the singularity theory of planar curves.  They also gave complete spelling rules for the RVT code in this case.

These studies were all concerned with the Monster tower whose base is $\mathbb R^2$.  In \cite{castro2}, we generalized this to towers with base $\mathbb R^n$.  We also discovered that this object was known to algebraic geometers as the Semple tower.  We also began the effort to generalize the RVT code, and find spelling rules to describe which words were admissible (Theorem \ref{Thm1}).  That effort is completed here (Theorem \ref{Thm2}).  The methods used in the present contribution were first developed in \cite{castro3}, in which we also classified points in the first four levels of the tower. Here, we will complete the spelling rules for base $\mathbb R^3$.  
Our techniques generalize to towers with base $\mathbb R^n$.

\subsection{Main Results}
The diffeomorphism group of $\mathbb R^3$ acts on the Monster tower, and the RVT code is an invariant labeling of orbits.  Note that the combinatorial data in the RVT code forces a finite number of inequivalent classes at each level of the tower, but there may be moduli within a given class (see \cite{MZ1}.
In \cite{castro2}, we stated the following incomplete spelling rules, which followed from \cite{MZ1}.

\begin{theorem}[\cite{castro2}]\label{Thm1}
In the Semple tower with base $\mathbb R^3$, every RVT code must begin with $R$, and $T_1$ cannot follow $R$.
\end{theorem}

Here, we add the missing rules, yielding the complete description of realizable RVT codes.  Our alphabet is the set $\{R, V, T_1, T_2, L_1, L_2, L_3\}$. Note that these seven letters correspond precisely to the seven possibilities found in Semple's original work \cite{semple}.  We therefore have the following combinatorial description of the diffeomorphism group orbits.

\begin{theorem}[Spelling Rules]
In the Semple tower with base $\mathbb R^3$, there exists at point $p$ with RVT code $\omega$ if and only if the word $\omega$ satisfies:
\begin{enumerate}\label{Thm2}
\item Every word must begin with $R$
\item $R$ must be followed by $R$ or $V$
\item $V$ and $T_1$ must be followed by $R, V, T_1$, or $L_1$
\item $T_2$ must be followed by  $R, V, T_2$, or $L_3$
\item $L_1, L_2$, and $L_3$ can be followed by any letter.
\end{enumerate}
\end{theorem}

For example, the word $RVVRVT_1L_1T_2L_3L_2$ is admissible, but $VT_2T_1RT_2$ breaks rules (1)--(4).
The following Table \ref{tab:rules} summarizes this Theorem.

\smallskip

\begin{table}
\caption{RVT Code Spelling Rules}
  \begin{tabular}{ |  l | c | r | }
    \hline
    Letter & Can be followed by & Cannot be followed by \\ \hline 
    $R$ & $R, V$ & $T_i, L_j$ \\ \hline 
    $V$ & $R, V, T_1, L_1$ & $T_2, L_2, L_3$ \\ \hline
    $T_1$ &$R, V, T_1, L_1$ & $T_2, L_2, L_3$ \\ \hline 
    $T_2$ & $R, V, T_2, L_3$ & $T_1, L_1, L_2$ \\ \hline 
    $L_1$ & $R, V, T_1, T_2, L_1, L_2, L_3$ & $\emptyset$ \\ \hline 
    $L_2$ & $R, V, T_1, T_2, L_1, L_2, L_3$ & $\emptyset$ \\ \hline 
    $L_3$ & $R, V, T_1, T_2, L_1, L_2, L_3$ & $\emptyset$ \\ \hline 
  \end{tabular}
  \label{tab:rules}
\end{table}

\smallskip

\subsection{Outline}
In Section 2, we give the requisite background material and references.  We define the Monster tower, Baby Monsters, and the RVT coding system.  

In Section 3, we describe our main tool, the method of critical hyperplanes.  We begin our main example which will inform the rest of the paper.  This example -- the code $RVL_1T_2$ -- will lend itself to a model proof of one spelling rule, whose technique can be repeated to obtain the remaining rules.  Moreover, this example will serve to demonstrate the ease with which our results could be extended to towers with bases $\mathbb R^n$ for $n>3$.  We choose this code to focus on because it neatly demonstrates the general method as well as some of the subtleties which abound in this work and thereby necessitate a delicate touch.  In particular, the code $RVL_1$ was studied extensively in \cite{castro3}, so we restate and build upon the work there.  We then amend the code by adding $T_2$, which is somewhat exotic and interesting but not overly complicated.


In Section 4, we restate our main theorem and attend to its proof.  We focus on one spelling rule, as the rest are proved in the same fashion, and the proofs are tedious.  
The main proof proceeds by induction on the number of letters appearing in the code which belong to the set $S=\{T_2, L_2, L_3\}$.  


\section{Background}

\subsection{The Tower}

The Monster/Semple tower is constructed through a series of {\it Cartan 
prolongations}.  Begin with a smooth $d$-dimensional manifold $M^0$ and a rank $r$ distribution (subbundle of $TM^0$) denoted $\Delta^0$.  The first prolongation is the fiber bundle
\[ M^1 = \bigcup_{p\in M^0}  \mathbb{P}\Delta^0_p,\]
whose elements have the form $(p,l)$, where $p$ is a point in $M^0$ and $l$ is a line in the subspace $\Delta^0_p$.  The distribution on $M^1$ is given by
\[ \Delta^1_{(p,l)}= (d\pi^1_0)^{-1}(l) \]
where $\pi^1_0 \colon M^1 \to M^0$ is the bundle projection.  Note that $M^1$ has dimension $d+r-1$, and that $\Delta^1$ is a rank $r$ distribution.

Iterating the prolongation procedure gives a sequence of manifolds
 \[ M^i = \bigcup_{p\in M^{i-1}}  \mathbb{P}\Delta_p^{i-1}.\]
  Every point in $M^i$ has the form $(p,l)$, where $p$ is a point in $M^{i-1}$ and $l$ is a line in the distribution $\Delta_p^{i-1}$. The dimension of $M^i$ is thus $d+i(r-1)$.  The bundle projection map $\pi^i_{i-1} \colon M^i \to M^{i-1}$ has fibers diffeomorphic to $\mathbb{P}\Delta_p^{i-1} \cong \mathbb{RP}^{r-1} $.  The rank $r$ distribution  on $M^i$ is given by
  \[ \Delta^i_{(p,l)}= (d\pi^i_{i-1})^{-1}(l). \]
  The distributions $\Delta^i$ are sometimes known as Goursat multi-flags.  
  
\begin{definition}  The \emph{Monster} or \emph{Semple tower} is the sequence of projective bundles
  \[ \dots \rightarrow M^i \rightarrow M^{i-1} \rightarrow \dots \rightarrow M^1 \rightarrow M^0 \]
   equipped with the distribution $\Delta^i$ at each level. 
\end{definition}

Of particular interest is the case of $M^0=\mathbb R^n$  and $\Delta^0=T\mathbb R^n$.
We refer to the consequent tower as the \emph{$\mathbb R^n$-tower} or the \emph{tower with $n$-dimensional base}.
The tower with base $M^0=\mathbb R^2$  and $\Delta^0=T\mathbb R^2$ has been studied extensively \cite{MZ1}.  Here, as in \cite{castro3}, we focus on the case $M^0=\mathbb R^3$  and $\Delta^0=T\mathbb R^3$.  However, our methods generalize to the $\mathbb R^n$-tower for arbitrary $n$.

To be clear, in the remainder of this paper we are taking $M^0=\mathbb R^3$  and $\Delta^0=T\mathbb R^3$.

\subsection{Regular, Critical, and Vertical Directions and Points}
  By composing the projection maps $\pi^k _{k-1},\  \pi^{k-1}_{k-2}, \ldots, \pi^{i+1}_{i}$ we obtain projections $\pi_i^k : M^k \to M^i$, $i < k$. For $p_k \in M^k$, we denote $\pi_i^k (p_k)$ by $p_i$.  The horizontal curves at level $i$  (tangent to $\Delta^i$)
naturally prolong (i.e., lift) to  horizontal curves at level $k$. 
However, the curves coinciding with fibers of $\pi_{i-1}^{i}$ are special -- they project down to points and 
are not prolongations of curves from below. They are called 
\emph{vertical} and can themselves be prolonged to (first order) 
tangency curves, then prolonged again to (second order)
tangency curves, and so on. Vertical curves and their prolongations are called 
{\it critical}. 
If a curve is vertical or critical then we say its tangent directions are as well.

Thus, at each level $i\geq 2$ there are vertical 
directions, and, in addition, at each level $i \geq 3$ there 
are tangency directions different from the vertical direction. At any level, 
all the remaining (non-critical) horizontal directions are called \emph{regular}. 
Finally, we call a point $(p,l)\in M^i$ regular, vertical, or critical if the direction of $l$ is.

\subsection{Baby Monsters and Critical Hyperplanes}

Recall that one can apply the prolongation procedure to any smooth manifold $F$ in place of $\mathbb R^3$.  In particular, we will prolong the fibers $F$ of the bundle projections $\pi_{i-1}^i$, obtaining new subtowers of the Monster tower.  We call these subtowers \emph{Baby Monsters}.  

Let $p_i \in M^i$ and consider the fiber $F_i(p_i):=(\pi_{i-1}^{i})^{-1}(p_{i-1})\subset M^i$.  This is an integral submanifold for $\Delta^i$, so we can prolong the pair $(F_i(p_i), TF_i(p_i))$.  Denote the $j$th prolongation of this pair by $(F^j_i(p_i), \delta_i^j)$.  Note that $F^j_i(p_i)$ is a smooth submanifold of $M^{i+j}$, and 
\[\delta_i^j(q)=\Delta^{i+j}(q) \cap T_qF^j_i(p_i)\]
for $q \in F^j_i(p_i)$. 

\begin{definition}
We call the tower $(F^j_i(p_i), \delta_i^j)$ a \emph{Baby Monster} born at level $i$.  For $q \in F^j_i(p_i)$, we call $\delta_i^j(q)$ a \emph{critical hyperplane}.
\end{definition}

Note that the Baby Monster is a subtower of the Monster tower, with $\dim F^j_i(p_i)=2+j $ and $\dim \delta_i^j(q)=2$.  While the terminology \emph{hyperplane} comes from a more general setting, here we will simply refer to \emph{critical planes}.

\subsection{KR Coordinates}
It is convenient to work in a canonical coordinate system, called \emph{Kumpera-Ruiz} or \emph{KR-coordinates} \cite{GKR}.  This is a generalization of jet coordinates for jet spaces, but that takes into account the projective nature of the fibers.  These coordinates were described in detail for the $\mathbb R^2$-tower in \cite{MZ1} and for our current case, the $\mathbb R^3$-tower, in \cite{castro2}.  We briefly summarize here for completeness, and refer the interested reader to Section 4.2 of \cite{castro2}.

The KR coordinates for $M^k$ are of the form $(x, y, z, u_1, v_1, \dots u_k, v_k)$.  They satisfy:
\begin{enumerate}
\item The projection $\pi^k_i(x, y, z, u_1, v_1, \dots u_k, v_k)=(x, y, z, u_1, v_1, \dots u_i, v_i)$;
\item The coordinates $u_k, v_k$ are affine coordinates for the fiber $F_k$;
\item There are $3^k$ many charts covering $M^k$, corresponding to the three affine charts needed to cover each $F_i \cong \mathbb{RP}^2$ for $1\leq i \leq k$.
\item 
The projective fiber $F_{i+1}$ is coordinatized homogeneously by $[df_i : du_i : dv_i]$, 
where $f_i$ is some coordinate from a lower level.  The covector $df_i$ is called the \emph{uniformizing coordinate} in \cite{castro3}.  Dividing the entries in $[df_{i}: du_{i} : dv_{i}]$ by one of the nonzero covectors yields local affine coordinates for the fiber $F_{i+1}$.  By convention, we always take $u_{i+1}$ to be the first (left-most) affine coordinate.    
\end{enumerate}

\

To illustrate $(4)$, consider the following two examples. If $[df_{i} : du_{i} : dv_{i}] = [1 : \frac{du_{i}}{df_{i}} : \frac{dv_{i}}{df_{i}}]$, then we take $u_{i+1} = \frac{du_{i}}{df_{i}}$ and $v_{i+1} = \frac{dv_{i}}{df_{i}}$.  Similarly, if $[df_{i} : du_{i} : dv_{i}] = [\frac{df_{i}}{du_{i}} : 1 : \frac{dv_{i}}{du_{i}}]$, then $u_{i+1} = \frac{df_{i}}{du_{i}}$ and $v_{i+1} = \frac{dv_{i}}{du_{i}}$. Detailed examples are worked below.




\subsection{RVT Codes}

We observed in \cite{castro2} that there are only three critical planes within each distribution $\Delta^i$.  The tangent space to the fiber is called the vertical plane; the other two arise as prolongations of vertical planes and are called tangency planes. In the most general setting, a tangency hyperplane is any hyperplane with nontrivial intersection with the vertical hyperplane.  In our setting, we have the following characterization.

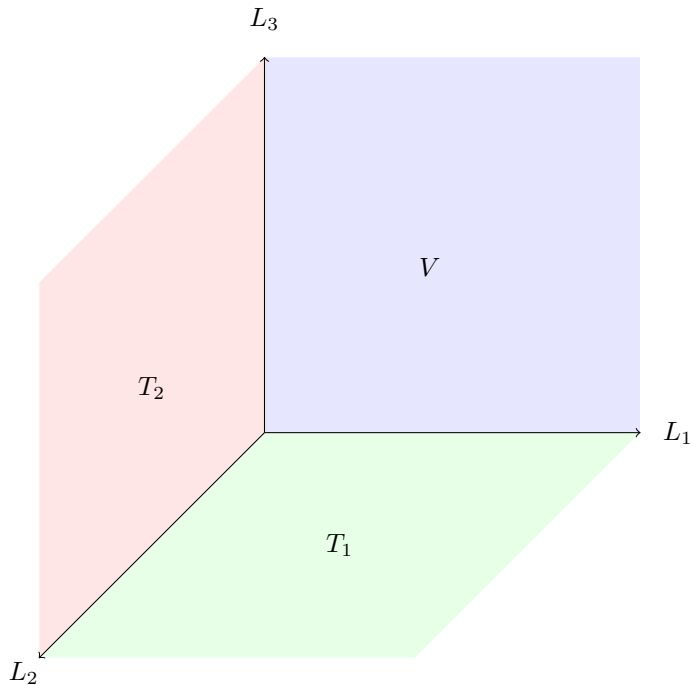
\begin{figure}

\begin{center}
\begin{tikzpicture}
\draw[fill=blue, draw=white, fill opacity=.1] (0,0) -- (5, 0) -- (5, 5) -- (0,5) -- (0, 0);
\draw[fill=green, draw=white, fill opacity=.1] (0,0) -- (-3, -3) -- (2, -3) -- (5,0) -- (0, 0);
\draw[fill=red, draw=white, fill opacity=.1] (0,0) -- (-3, -3) -- (-3, 2) -- (0,5) -- (0, 0);
\draw[->] (0,0) -- (5,0);
\draw[->] (0,0) -- (0,5);
\draw[->] (0,0) -- (-3,-3);
\node at (5.5,0) {$L_1$};
\node at (-3.2,-3.2) {$L_2$};
\node at (0,5.5) {$L_3$};

\node at (1,-1.5) {$T_1$};
\node at (-1.5,.6) {$T_2$};
\node at (2.2,2.2) {$V$};
\end{tikzpicture}
\end{center}
\caption{The three critical planes $V, T_1$, and $T_2$, and their intersections, the distinguished lines $L_1, L_2$, and $L_3$.}
\label{planes}
\end{figure}

\begin{definition}\label{letters} Let $q\in M^i$.

\begin{enumerate}

\item The \emph{vertical plane} $V(q)$ is the critical plane $\delta_i^0(q)=T_qF_i(q)$.  In KR-coordinates, $V(q)=\sspan\{\frac{\partial}{\partial u_i},\frac{\partial}{\partial v_i}\}$.  It is given projectively by $[df_i : du_i : dv_i]=[0 : a : b]$ for $a, b \in \mathbb R$.

\item The plane $T_1(q)$ is the unique critical plane in $\Delta^i$ which intersects $\sspan\{\frac{\partial}{\partial v_i}\}$.  It is given projectively by $[df_i : du_i : dv_i]=[a : 0 : b]$ for $a, b \in \mathbb R$.

\item The plane $T_2(q)$ is the unique critical plane in $\Delta^i$ which does not intersect $\sspan\{\frac{\partial}{\partial v_i}\}$.  It is given projectively by $[df_i : du_i : dv_i]=[a : b : 0]$ for $a, b \in \mathbb R$.

\item The distinguished lines $L_j(q)$ for $j=1,2,3$ are given by:
\begin{enumerate}[(i)]
\item $L_1=V \cap T_1$
\item $L_2=T_1 \cap T_2$
\item $L_3=V \cap T_2$
\end{enumerate}
\end{enumerate}
See Figure \ref{planes}.
\end{definition}

In this definition, we often drop the explicit dependence on $q$ when the context is clear.  Also, in homogeneous coordinates, we cannot have $a$ and $b$ both zero, and we will usually assume without loss of generality that $a\neq 0$.  Finally, we clarify the terminology.  Here $V(q)$ is a linear subspace of $\Delta^i(q) \subset T_qM^i$.  When working in homogeneous coordinates, we are identifying this plane with $\mathbb P V(q)\subset \mathbb P \Delta^i(q) \subset M^{i+1}$.  Similarly for the other planes and lines in this definition.  Note again that this definition has analogue in \cite{semple}.

Now a point $p_{i+1}=(p_i, l_i)$ is assigned a letter from $\{R, V, T_1, T_2, L_1, L_2, L_3\}$ according to whether $l_i$ lies in one of the critical planes or distinguished lines given in Definition \ref{letters}.  Here, the lines $L_j$ take precedence, so $l_i$ lying in $L_3$ is assigned the letter $L_3$, even though it also lies in both $V$ and $T_2$.  If $l_i$ does not lie in any of these, then it is regular (see above) and assigned the letter $R$. If $l_i$ is assigned the letter $\alpha$, then we say that $p_{i+1}$ is an \emph{$\alpha$ point.}
Note that in \cite{castro2}, the letters $T_2, L_2$, and $L_3$ were unknown, and the notation was $T=T_1$ and $L=L_1$.  All letters besides $R$ are called \emph{critical letters}.

\begin{definition}
The \emph{RVT code} of a point $p\in M^k$ is a word $\omega=\omega_1\omega_2\dots \omega_k$ in the letters $\{R, V, T_1, T_2, L_1, L_2, L_3\}$, where $\omega_i=\alpha$ if $\pi_i^k(p)$ is an $\alpha$ point.

\end{definition}

\begin{example}\label{ex1}
Suppose $p_3 \in M^3$ has RVT code $\omega=RVL_1$.  This means that $p_3=(p_2, l_2)$ with $l_2 = L_1(p_2)$, and $p_2=(p_1, l_1)$ with $l_1 \subset V(p_1)$.  Every direction in $\Delta^1$ is regular, so the leading letter $R$ yields no information.
\end{example} 

For convenience, sometimes we will also denote by $\omega$ the set of all points with RVT code $\omega$.  For example, we may write $p\in RVL_1T_2$ to signify that $p$ has RVT code  $RVL_1T_2$.

This coding provides a coarse stratification of points in the Monster/Semple tower.  Recall that finite jets of diffeomorphisms act on the tower.  Points which lie in the same orbit must have the same RVT code.  However, there may exist multiple orbits within the same RVT strata.  For details, see \cite{castro3} or \cite{morpel}. 

\section{The Critical Hyperplane Method}

\subsection{Configurations}

\begin{figure}
\centering
\begin{tabular}{cc}
\includegraphics[width=0.177\textwidth]{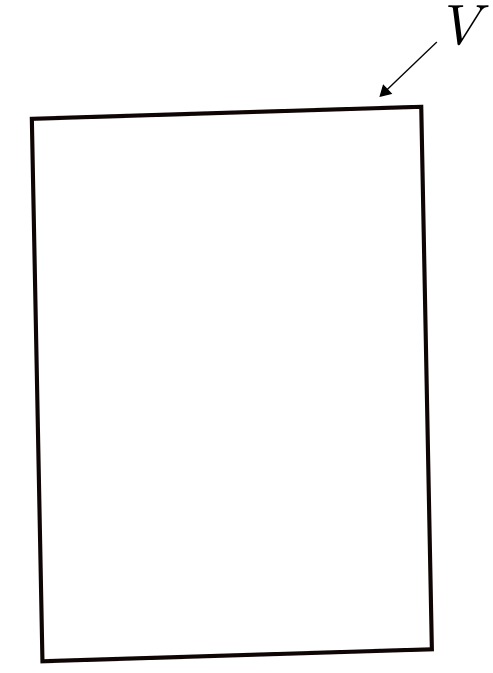} \hspace{1in}
\includegraphics[width=0.25\textwidth]{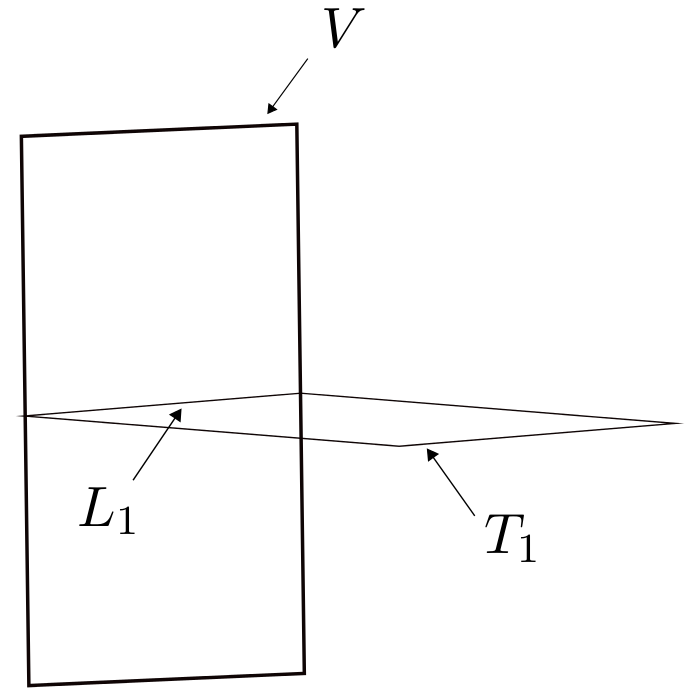} \\ \\
\includegraphics[width=0.24\textwidth]{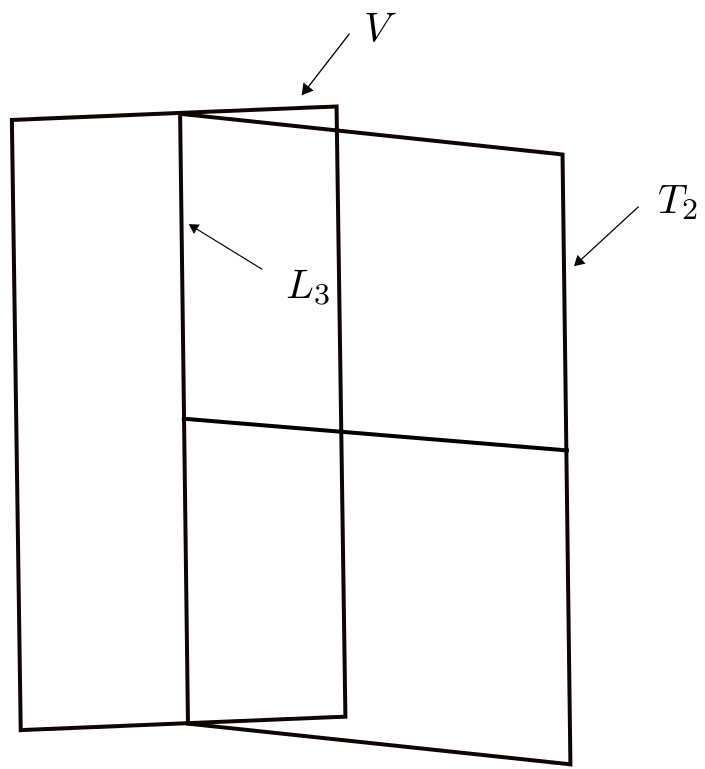} \hspace{.5in}
\includegraphics[width=0.25\textwidth]{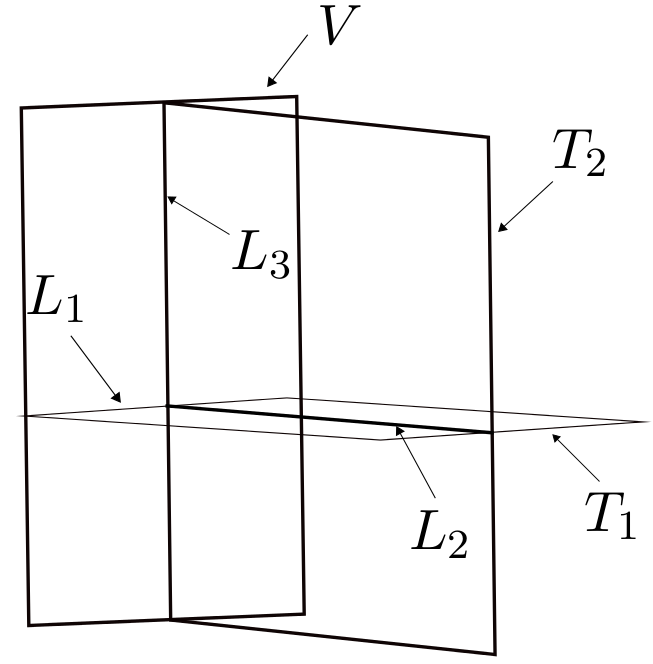} 
\end{tabular}
\caption{Critical plane configurations that can appear in the distribution above an $R$ point (top left), a $V$ or $T_1$ point (top right), a $T_2$ point (bottom left), and an $L_j$ point (bottom right).}
\label{fig1}
\end{figure}

This method relies on the non-trivial fact that certain critical planes appear over certain points, while others may not.  In particular, there are four possible configurations over a point $p\in M^k$; these are shown in Figure \ref{fig1}.  We will show how each configuration is possible only if $p$ belongs to certain RVT classes.    Specifically, we have Table \ref{tab:critplanes}, which is effectively equivalent to Theorem \ref{Thm2}.  
Note that saying that $p$ is an $\alpha$ point is the same as saying that $\alpha$ is the last letter in the RVT code for $p$.


\begin{table}
  \caption{Critical Hyperplane Configurations}
  \begin{tabular}{ |  c | c | }
    \hline
    Last letter in RVT code of $p\in M^k$ & Critical planes appearing in $\Delta^k(p)$\\ \hline 
    $R$ & $V$\\ \hline 
     $V$ or $T_1$ & $V$ and $T_1$\\ \hline 
      $T_2$ & $V$ and $T_2$\\ \hline 
       $L_1, L_2,$ or $L_3$ & $V, T_1$, and $T_2$\\ \hline 
  \end{tabular}  
  \label{tab:critplanes}  
\end{table}


The remainder of this paper will be dedicated to explaining why these possibilities are exhaustive.

\subsection{The Method}

We now describe the explicit method from which we derive all our results.  This will be applied to specific examples shortly.  The critical hyperplane method was implicit in parts of \cite{MZ1}, made explicit in \cite{castro2}, exploited for the classification problem in \cite{castro3}, and is perfected here.  This gives a blueprint for characterizing all Baby Monsters and determining all spelling rules for the $\mathbb R^n$ tower for any $n$.  

Begin with an RVT code $\omega$ of a point $p\in M^k$.  We wish to understand which critical letters can be added to the end of the code (one can always trivially add the letter $R$).  In order to do so, we must understand which critical planes lie above $p$.  Since critical planes live within Baby Monsters, we must determine which Baby Monsters are present, and for those which are, we seek to find the levels at which they were born.  

We first determine the local KR-coordinate chart containing $p$.  We can then describe the distribution $\Delta^k(p)$ in coordinates.  We then choose a critical plane $V, T_1$, or $T_2$, write it in coordinates as in Definition \ref{letters}, and trace the coordinate representations backwards, projecting down to lower levels of the tower, one at a time. 

If at some level $i$ we find that both fiber coordinates $u_i$ and $v_i$ are non-vanishing, then our critical plane must arise as the prolongation of the vertical plane $V_i$.  Our critical plane therefore lives in the Baby Monster born at level $i$, and is equal to $\delta_i^k(p)$.  This would confirm that the critical plane we chose indeed appears in $\Delta^k(p)$.

If, however, we reach the base without finding such a Baby Monster, then the plane we chose cannot exist in $\Delta^k(p)$.  We can shorten the procedure of tracing each plane back to the base by using previously established configuration possibilities and proceeding inductively.  

While this is not an algorithm in the strictest sense, it can theoretically determine which configurations are possible above any given point.  As one might suspect, this can at times become extremely tedious, and would not be particularly enlightening for the reader.  For this reason, we will focus the remainder of the paper on a few specific examples to demonstrate the efficacy of the method for determining spelling rules, while skipping some of the routine verification that was required to complete our results. 

It is obvious that the vertical plane $V$ appears above every point -- it is just the tangent space to the fiber.  So in the method just described, we need only focus on whether or not $T_1$ and $T_2$ exist (here, since we are concerned with the $\mathbb R^3$-tower -- one immediately sees how this method generalizes to the $\mathbb R^n$-tower).  Some of our results here (those needed for the proof of Theorem \ref{Thm2}) are summarized near the end of the paper in Table \ref{tab:codes}.  We will prove some of these relations here -- the rest are obtained by identical methods.

\begin{example}[$RVL_1$]\label{ex: RVL1}
We continue investigating the case begun in Example \ref{ex1}.  Suppose $p_3 \in M^3$ has RVT code $\omega=RVL_1$.


\begin{table}
\caption{Summary of Example \ref{ex: RVL1}: $RVL_1$}
  \begin{tabular}{ |  p{1.1cm} | p{4.1cm} | c | p{3.1cm} | p{3.3cm} | }
    \hline
    Level $i$ & Coordinates on $M^i$ &  $\mathbb P\Delta^{i-1} = F_i$ coords. & Critical planes in $\Delta^i$ & RVT code of $p_i$ \\ \hhline{|=|=|=|=|=|}
    
$0$ & $(x, y, z)$ & n/a & none & n/a \\ \hline 

$1$ & $(x, y, z, u_1, v_1) \newline u_1=\frac{dy}{dx}, v_1=\frac{dz}{dx}$ & $[dx : dy : dz]$ & $V(p_1)=\delta_1^0$ & $p_1=(p_0, l_0)\in R \newline l_0 \subset \Delta^0=T_{p_0}M^0$ \\ \hline 

$2$ & $(x, y, z, u_1, v_1, u_2, v_2) \newline u_2=\frac{dx}{du_1}, v_2=\frac{dv_1}{du_1}$ & $[dx : du_1 : dv_1]$ & $V(p_2)=\delta_2^0, \newline T_1(p_2)=\delta_1^1$ &  $p_2=(p_1, l_1)\in RV \newline l_1 \subset V(p_1) \subset \Delta^1$ \\ \hline 

$3$ & $(x, y, z, u_1, v_1, u_2, v_3, u_3, v_3) \newline u_3=\frac{du_1}{dv_2}, v_3=\frac{du_2}{dv_2}$ & $[du_1 : du_2 : dv_2]$ & $V(p_3)=\delta_3^0, \newline T_1(p_3)=\delta_2^1, \newline T_2(p_3)=\delta_1^2$ & $p_3=(p_2, l_2)\in RVL_1 \newline l_2 = L_1(p_2) \subset \Delta^2$ \\ \hline 
  
  \end{tabular}
  \label{tab:ex2}
\end{table}


\begin{figure}
\centering
\includegraphics[width=0.45\textwidth]{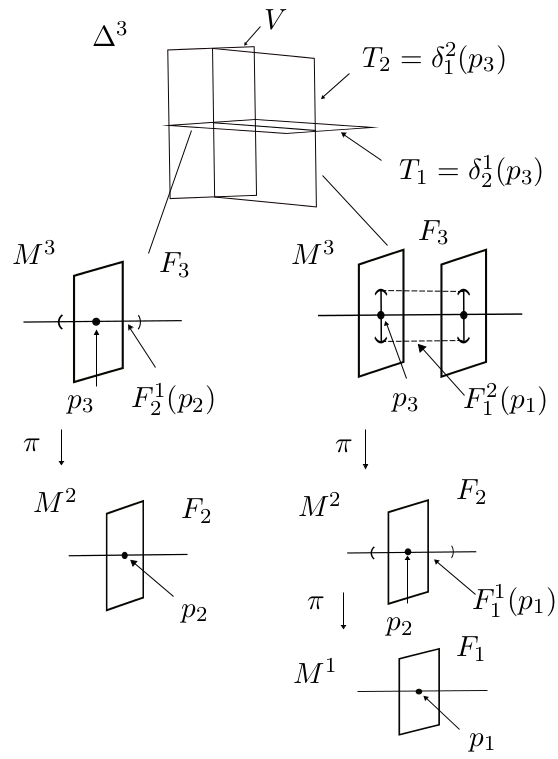}
\caption{Critical plane configuration over $p_3\in RVL_1$. The left side shows the birth of $T_1(p_3)=\delta_2^1(p_3)$ as the first prolongation of the vertical plane at level 2. The right side shows the birth of $T_2(p_3)=\delta_1^2(p_3)$ as the second prolongation of the vertical plane at level 1.  These two Baby Monsters meet in $\Delta^3$, and their intersection is the distinguished line $L_2(p_3)$.  See Example \ref{ex: RVL1}.}  
\label{fig2}
\end{figure}

\subsubsection*{Level 1.}
Begin with the global coframe $\{dx, dy, dz\}$ for $\Delta^0=T\mathbb R^3$.  Our chart will be centered at $p_1=(p_0, l_0) \in M^1$ where $p_0=(0,0,0)$.  Introduce affine fiber coordinates $[dx : dy : dz]$ on $F_1(p_1)$.   Without loss of generality, assume $dx|_{l_0} \neq 0$. Then $[dx : dy : dz]=[1: \frac{dy}{dx} : \frac{dz}{dx}]$.  Now let 
\[u_1= \frac{dy}{dx}, \qquad v_1= \frac{dz}{dx}\]
so that $$\Delta^1(p_1)=\{dy-u_1dx=0,\ dz-v_1dx=0 \}.$$

\subsubsection*{Level 2.} 
Since $l_1\subset V(p_1)=\sspan\{\frac{\partial}{\partial u_1}, \frac{\partial}{\partial v_1}\}$, we know $l_1= \sspan\{a\frac{\partial}{\partial u_1}+b\frac{\partial}{\partial v_1}\}$ with $a,b$ not both zero.  Without loss of generality, assume $a\neq0$.  Then near this point we have $[dx : du_1 : dv_1]=[\frac{dx}{du_1}: 1 : \frac{dv_1}{du_1}]$.  This yields the affine coordinates 
\[u_2= \frac{dx}{du_1}, \qquad v_2= \frac{dv_1}{du_1}\]
so that $$\Delta^2(p_2)=\{dy-u_1dx=0,\ dz-v_1dx=0,\ dx-u_2du_1=0,\ dv_1-v_2du_2=0 \}.$$

\subsubsection*{Level 3.} 
Now $l_2=L_1(p_3)=V(p_2) \cap T_1(p_3)$, so we want coordinate representations of the $V$ and $T_1$ planes in $\Delta^2(p_2)$.  According to Definition \ref{letters}, $V(p_2)$ is given by $du_1=0$ and $T_1(p_2)$ is given by $du_2=0$, so we have $du_1|_{l_2}=0$ and $du_2|_{l_2}=0$. This forces our coordinates near $p_3$ to have the form 
 $[du_1 : du_2 : dv_2]=[\frac{du_1}{dv_2} : \frac{du_2}{dv_2}: 1]$.  This yields the affine coordinates 
\[u_3= \frac{du_1}{dv_2}, \qquad v_3= \frac{du_2}{dv_2}\]
so that 

\begin{align*}
\Delta^3(p_3)= \{& dy-u_1dx=0,\ dz-v_1dx=0,\  dx-u_2du_1=0,\\  
& dv_1-v_2du_2=0,\ du_1-u_3dv_2=0,\ du_2-v_3dv_2=0  \}.
\end{align*}

This completes the first step of the process, as we have determined the local KR-coordinates around $p_3$ and described the distribution $\Delta^3(p_3)$ in these coordinates.  Note that here $p_3=(0,0,0,0,0,0,0,0,0)$.

\subsubsection*{Appearance of $T_1$.}  
We now determine which of the critical planes $T_1$ and $T_2$ lie above $p_3$ in $\Delta^3(p_3)$, which is coframed\footnote{Technically, this coframes the projectivized space.  But as we often identify $\Delta^k(p)\subseteq T_{p} M^k$ with $\mathbb P \Delta^k(p)\subseteq M^{k+1}$, this abuse of notation is convenient and should not cause confusion.} by  $[dv_2 : du_3 : dv_3]$.  
First consider $T_1$, given by  $[a : 0 : b]$ with $a\neq 0$.  We assume for now that it exists within some Baby Monster, and we will either find this Baby Monster or derive a contradiction.  Since $[dv_2 : du_3 : dv_3]=[a : 0 : b]$ here with $a\neq 0$, we see that $u_3$ is identically zero on the Baby Monster, while $v_2$ and $v_3$ are not.  Now, since $\Delta^2$ is coframed by $[du_1 : du_2 : dv_2]$ near $p_2$, and since $u_3= \frac{du_1}{dv_2}$ and $v_3= \frac{du_2}{dv_2}$, this forces the Baby Monster to have the form
$[du_1 : du_2 : dv_2]=[0 : c : d]$.  Since this is the form of a vertical plane, we can stop and conclude that $T_1(p_3)$ exists, and lies inside the Baby Monster born at level 2.  That is, the plane $T_1(p_3)=\delta_2^1$, which is the first prolongation of the tangent space to the fiber $F_2(p_2)$.

\subsubsection*{Appearance of $T_2$.} 
Next, we repeat this process for $T_2$, given by  $[a : b : 0]$ with $a\neq 0$.  We assume for now that it exists within some Baby Monster, and we will either find this Baby Monster or derive a contradiction.  Since $[dv_2 : du_3 : dv_3]=[a : b : 0]$ here with $a\neq 0$, we see that $v_3$ is identically zero on the Baby Monster, while $v_2$ and $u_3$ are not.  Now, since $\Delta^2$ is coframed by $[du_1 : du_2 : dv_2]$ near $p_2$, and since $u_3= \frac{du_1}{dv_2}$ and $v_3= \frac{du_2}{dv_2}$, this forces the Baby Monster to have the form $[du_1 : du_2 : dv_2]=[c : 0 : d]$.  Note that unlike the previous case, this is not vertical, so we must continue searching another level down.
Since $\Delta^1$ is coframed by $[dx : du_1 : dv_1]$ near $p_1$, and since $u_2= \frac{dx}{du_1}$ and $v_2= \frac{dv_1}{du_1}$, this forces the Baby Monster to have the form $[dx : du_1 : dv_1]=[0 : e : f]$.
Since this is the form of a vertical plane, we can stop and conclude that $T_2(p_3)$ exists, and lies inside the Baby Monster born at level 1.  That is, the plane $T_2(p_3)=\delta_1^2$, which is the second prolongation of the tangent space to the fiber $F_1(p_1)$.

\subsubsection*{Summary}
We conclude that both planes $T_1$ and $T_2$ occur above a point with RVT code $\omega=RVL_1$, so that both codes $RVL_1T_1$ and $RVL_1T_2$ are admissible and realized (assuming temporarily that $\omega$ is admissible).  Compare this result with Theorem \ref{Thm2} and Figure \ref{fig1}. 
Also see Figure \ref{fig2} for an illustration of this situation.  We summarize the results of this example in Table \ref{tab:ex2}.

\end{example}

\begin{figure}
\centering
\includegraphics[width=.4\textwidth]{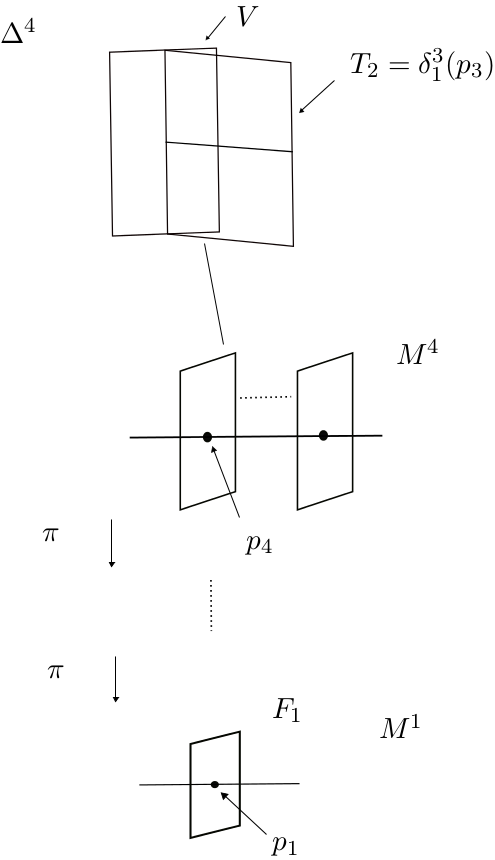}
\caption{Critical plane configuration over $p_4\in RVL_1T_2$. This shows the birth of $T_2(p_4)=\delta_1^3(p_4)$ as the third prolongation of the vertical plane at level 1.  See Example \ref{ex: RVL1T2}.}
\label{fig3}
\end{figure}

\begin{example}[$RVL_1T_2$]\label{ex: RVL1T2}
We continue the work from the previous example, and consider the case of $p_4\in M^4$ with RVT code $RVL_1T_2$.  This is admissible by the preceding computations, and indeed, all results from that example hold here.  As the general techniques were made explicit there, we omit some tiresome details here.

First, one finds affine coordinates $u_4= \frac{du_3}{dv_2}$ and $v_4= \frac{dv_3}{dv_2}$ for the fiber $F_4(p_4)$.  Next, recall that  $\Delta^3(p_3)$, is coframed by  $[dv_2 : du_3 : dv_3]$, and $T_2(p_3)$ locally satisfies $dv_3=0$, with $dv_2$ non-vanishing and $du_3$ not identically zero.  This implies that $v_4=0$, but $u_4$ is non-zero. (If $u_4(p_4)$ were zero, then there would be no vertical component, and $l_3$ would lie in a regular direction instead of in $T_2$.)

Second, we show that $T_2$ does occur in $\Delta^4(p_4)$. This computation is nearly identical to those presented in the previous example, so we omit it.  One finds that $T_2(p_4)=\delta_1^3$.

Finally, we show that $T_1$ cannot occur in $\Delta^4(p_4)$.  
If it did, it would have the form $[dv_2 : du_4 : dv_4]=[a : 0 : b]$ with $a \neq 0$. But $p_4=(p_3,l_3)$ with $l_3\subset \delta_1^2 = T_2(p_3)$.   This implies $du_4|_{l_3}=0$, so $u_4(p_4)=0$, which contradicts the fact that $u_4$ is non-zero in a neighborhood of $p_4.$  

We have shown that the $T_2$ critical plane occurs, but $T_1$ does not, in $\Delta^4(p_4)$ for $p_4$ in the class $RVL_1T_2$.  We conclude that the code $RVL_1T_2$ can be amended with letters $R, V, T_2$, and $L_3$, but not with $T_1, L_1,$ or $L_2$.  Compare with Theorem \ref{Thm2}, Figure \ref{fig1}, and the second row of Table \ref{tab:codes}. 
Also see Figure \ref{fig3} for an illustration of this situation.

\end{example}


\section{Spelling Rules}

In this section we will outline the proof of Theorem $\ref{Thm2}$ from the Introduction, which we restate here.

\newtheorem*{Thm2}{Theorem \ref{Thm2}}
\begin{Thm2}[Spelling Rules]
In the Semple tower with base $\mathbb R^3$, there exists at point $p$ with RVT code $\omega$ if and only if the word $\omega$ satisfies:
\begin{enumerate}
\item Every word must begin with $R$
\item $R$ must be followed by $R$ or $V$
\item $V$ and $T_1$ must be followed by $R, V, T_1$, or $L_1$
\item $T_2$ must be followed by  $R, V, T_2$, or $L_3$
\item $L_1, L_2$, and $L_3$ can be followed by any letter.
\end{enumerate}
\end{Thm2}

Let us begin with an overview of the method of proof.  The first two rules are well known and appear in \cite{castro2} and \cite{castrothesis}.  Rule (3) can be checked by direct calculation; this is tedious but straightforward and we omit the computation here.  The same can be said for the part of rule (5) concerning the letter $L_1$.  The technique is illustrated by examples in \cite{castro3} and the three examples above.  For example, one finds that for any point $p \in \lambda L_1$, the plane $T_1(p)$ is obtained by prolonging the vertical plane from one level below.  In other words, $T_1(p)=\delta_{k-1}^1$.  Similarly, the plane $T_2(p)$ is the prolongation of the $T_1$ plane from one level below.  This is independent of the code $\lambda$.

To prove the remaining rules, (4) and most of (5), 
we proceed by induction on the number of letters $T_2, L_2$, or $L_3$ appearing in the code.  This proof is more delicate.  Set $S=\{T_2, L_2, L_3\}$.  For the base case, we must prove that the spelling rules hold for an RVT code $\omega$ containing only one letter $\alpha \in S$.  For the inductive step, we must prove that the spelling rules hold for an arbitrary code $\omega$, using the inductive hypothesis that the rules hold for any code containing fewer letters $\alpha \in S$.  In both steps, we assume without loss of generality that the letter $\alpha$ appears at the end of the code in question.

Unfortunately (but perhaps unsurprisingly given the examples above), this method requires investigating a large number of specific cases, as well as a considerable number of tedious calculations.  We therefore spare the reader details of all cases, and the lengthy but routine computations which are required to prove each spelling rule rigorously.  Instead, we will focus in detail on one particular rule: the fourth.  
We hope that this approach will yield sufficient detail to introduce the mechanics of the method to the reader, while sparing the reader dozens of pages of nearly identical calculations.  We chose these particular cases as they exhibit generally typical behavior, but with a few of the subtleties which necessitate special care and patience.

\subsection{Base Case}
We assume rules (1) -- (3) have been proved.  
Here we will provide details for rule (4); the remaining proofs are very similar.
To this end, let $\omega$ be an RVT code of length $k$, ending with the letter $T_2$.  We will show that codes $\omega R, \omega V, \omega T_2$, and $\omega L_3$ do occur at level $k+1$, while $\omega T_1, \omega L_1$, and $\omega L_2$ are impossible.  We prove this by induction on the number of letters $\alpha \in S=\{T_2, L_2, L_3\}$ appearing in $\omega$.

We first prove the base case.  Assume $\omega=\lambda T_2$, where $\lambda$ does not contain any letter from $S$.  We will show that rule (4) holds for this $\omega$.  We prove this by considering the potential letters preceding $T_2$.  By rules (2) and (3), $T_2$ cannot be preceded by $R$  or $V$ or $T_1$.  Since we have assumed that $\lambda$ contains no letters from $S$, we know $T_2$ cannot be preceded by $T_2, L_2$, or $L_3$.  We therefore consider the only remaining possibility: $T_2$ is preceded by $L_1$.  Note that for convenience we will use $\lambda$ to denote any sub-code of $\omega$, regardless of its length.


So we proceed assuming our code has the form $\omega=\lambda L_1T_2$, where $\lambda$ contains no elements from $S$.  
Thus, the predecessor of $L_1$ can only be $V, T_1$, or $L_1$. We have three possible cases.

\subsubsection*{Case 1: $\omega=\lambda VT_1^mL_1T_2, m\geq0$.}
Assume our code has length $k$ and is of the form $\omega=\lambda VT_1^mL_1T_2$ with $m\geq0$.  If $m=0$, then $V$ precedes $L_1$; if $m>1$, then $T_1$ does.  The third possibility, where $L_1$ precedes $L_1$, is treated as a separate case below.

In fact, we can assume without loss of generality that $\omega=RVT_1^mL_1T_2$.  This is valid because the plane $T_2(p_k)$ is the (possibly multi-step) prolongation of some vertical plane from a lower level.  That is, $T_2(p_k)=\delta_i^j$ for some Baby Monster, and this subtower could not have been born at a level below the last letter $V$ in the RVT code. 

Now consider $\omega=RVT_1^mL_1T_2$.  We have $k=m+4$.  We wish to show that the spelling rules hold for $\omega$.  This is to show that
the codes $\omega \alpha$ are realized for $\alpha=R, V, T_2, L_3$, but are impossible for $\alpha = T_1, L_1, L_2$.  Now there are regular and vertical directions in each distribution plane, so it is clear that $\alpha=R$ or $V$ are possible.   Recall from Definition \ref{letters} that $L_1=V \cap T_1$, $L_2=T_1 \cap T_2$, and $L_3=V \cap T_2$.  It is therefore sufficient to simply show that $\alpha=T_2$ is possible, while $\alpha=T_1$ is not.  

The proof here is nearly identical to that provided in Example \ref{ex: RVL1T2}.  In fact, that example gives precisely the case where $m=0$.  Recall that in that case, $T_1$ could not appear and  $T_2(p_4)=\delta_1^{3}$.  For $m>1$, we easily verify that, again $T_1$ cannot appear, and  $T_2(p_{m+4})=\delta_1^{m+3}$.  The key observation is the following. The vertical plane $V(p_1)$ is coframed by $[dx : du_1 : dv_1]=[a: b: 0]$ with $a \neq 0$.   The prolongation of this plane is $T_1(p_2)=\delta_1^1$, which is coframed by $[du_1:du_2:dv_2]=[a:0:b]$ with $a\neq0$.  For $m>1$, we continue this process and find that the $m$th prolongation of $V(p_1)$ is $T_1(p_{m+1})=\delta_1^{m}$, which is coframed by $[du_1:du_{m+1}:dv_{m+1}]=[a:0:b]$ with $a\neq0$.  The rest of the steps are the same as in Example \ref{ex: RVL1T2}.

\subsubsection*{Case 2: $\omega=\lambda L_1T_1^mL_1T_2, m\geq1$.}
This case is nearly identical to the previous.  Here, one finds again that the vertical plane in $\Delta^{k-m-3}$ prolongs $m+3$ times to give the plane $T_2(p_k)$.

\subsubsection*{Case 3: $\omega=\lambda L_1L_1T_2$.}
The method here is the same as in Case 1, so we will omit some of the readily checked details.  Again suppose the length of $\omega$ is $k$.  Then $\Delta^k$ is coframed by $[df_k : du_k : dv_k]$, and $T_2(p_k)$ would have the form $[df_k : du_k : dv_k]=[a:b:0]$ with $a \neq0$ and $df_k=dv_{k-2}$.  Its projection in $\Delta^{k-1}$ will have the form $[df_{k-1} : du_{k-1} : dv_{k-1}]=[a:b:0]$ with $a \neq0$ and $df_{k-1}=dv_{k-2}$.  Its projection in $\Delta^{k-2}$ will have the form $[df_{k-2} : du_{k-2} : dv_{k-2}]=[a:0:b]$ with $a \neq0$ and $df_{k-2}=dv_{k-3}$.  Finally, its projection in $\Delta^{k-3}$ will have the form $[df_{k-3} : du_{k-3} : dv_{k-3}]=[0:a:b]$ with $a \neq0$.  At this point, we can see that this is the vertical plane $V(p_{k-3})$, so we find that $T_2(p_k)$ does indeed exist in $\Delta^k$, and that it is equal to $\delta_{k-3}^3$.

A computation similar to this one and those found in Example \ref{ex: RVL1T2} shows that $T_1(p_k)$ cannot exist.  In short, one repeats this computation beginning with $T_1(p_k)$ of the form $[df_k : du_k : dv_k]=[a:0:b]$ with $a \neq0$, and at some point a contradiction is obtained in that some coordinate is forced to be both zero and nonzero. 


\medskip

This establishes the base case for the proof of rule (4) by induction.  We showed that rule (4) holds for any RVT code containing a single member of $S$ (which, in the context of rule (4), must naturally be the letter $T_2$.)  
These three cases comprise the top three rows in Table \ref{tab:codes}.  The remaining cases are displayed as the lower six rows in Table \ref{tab:codes}; their proofs are similar.

\begin{table}
\caption{Base Cases of Inductive Proof}
\begin{tabular}{|c|c|c|}
\hline
RVT code of $p_k\in M^k$ &      $T_{1}(p_k)$ & $T_{2}(p_k)$ \\
\hhline{|=|=|=|}







  $ \lambda V T_1^{m} L_1 T_{2} \, \, \text{for} \, \,  m \geq 0$  &          None                                    & $\delta ^{m+3}_{k-m-3}(p_{k })$      \\
    $ \lambda L_1 T_1^{m} L_1 T_{2} \, \, \text{for} \, \,  m \geq 1$  &          None                                    & $\delta ^{m+3}_{k-m-3}(p_{k })$      \\
     $ \lambda L_1 L_1 T_{2}$                                                           &          None                                    & $\delta ^{3}_{k-3}(p_{k})$                  \\
\hhline{|=|=|=|}

  
    $ \lambda VT_1^{m} L_1L_{2} \, \, \text{for} \, \, m \geq 0$           &        $\delta ^{2}_{k-2}(p_{k})$      &  $\delta ^{m+3}_{k-m-3}(p_{k})$             \\
     $ \lambda L_1T_1^{m} L_1L_{2} \, \, \text{for} \, \, m \geq 1$           &        $\delta ^{2}_{k-2}(p_{k})$      &  $\delta ^{m+3}_{k-m-3}(p_{k})$             \\
      $ \lambda L_1L_1L_{2}$                                                                 &       $\delta ^{2}_{k-2} (p_{k})$      & $\delta ^{3}_{k-3}(p_{k})$            \\
\hhline{|=|=|=|}

  
     $ \lambda VT_1^{m} L_1 L_{3} \, \, \text{for} \, \, m \geq 0$           &        $\delta ^{1}_{k-1}(p_{k})$ &   $\delta ^{m + 3}_{k-m-3}(p_{k})$  \\
      $ \lambda L_1T_1^{m} L_1 L_{3} \, \, \text{for} \, \, m \geq 1$           &        $\delta ^{1}_{k-1}(p_{k})$ &   $\delta ^{m + 3}_{k-m-3}(p_{k})$  \\
       $ \lambda L_1L_1 L_{3}$                                                                &        $\delta ^{1}_{k-1}(p_{k})$ &   $\delta ^{3}_{k-3}(p_{k})$   \\
\hline

\end{tabular}
\label{tab:codes}
\end{table}

\subsection{Inductive Step}

We now take $\omega$ to be an arbitrary RVT code of length $k$.  We assume that $\omega$ ends with some letter from $S$, and we will show that the spelling rules hold for $\omega$.  Our inductive hypothesis states that the spelling rules hold for any code which contains fewer letters from $S$ than $\omega$ does.  

As above, we will focus on rule (4), so our code should end with the letter $T_2$. So we have $\omega=\lambda T_2$ and our inductive hypothesis allows the assumption that $\lambda$ satisfies the spelling rules.  
We wish to show that, at level $k+1$, the codes $\omega \alpha$ are realized for $\alpha=R, V, T_2, L_3$, but are impossible for $\alpha = T_1, L_1, L_2$.  Now there are regular and vertical directions in each distribution plane, so it is clear that $\alpha=R$ or $V$ are possible.   Recall from Definition \ref{letters} that $L_1=V \cap T_1$, $L_2=T_1 \cap T_2$, and $L_3=V \cap T_2$.  It is therefore sufficient to simply show that $\alpha=T_2$ is possible, while $\alpha=T_1$ is not.  

Now since $\lambda$ clearly has (exactly one) fewer letters from $S$ than $\omega$ does, it must obey the spelling rules by assumption.  So $T_2$ must be preceded by either $T_2, L_1, L_2$, or $L_3$.  There are four cases here, but we will give details for just the first and second.  The other two are nearly identical. 

\subsubsection*{Case 1: $\omega=\lambda T_2T_2$}
Suppose $p_k\in M^k$ has RVT code $\omega=\lambda T_2T_2$.
From the discussion above, it suffices to prove that $T_2$ appears in $\Delta^k$, while $T_1$ does not.  Now the distribution $\Delta^{k}$ is coframed by $[df_{k} : du_{k} : dv_{k}]$.  Two levels down,  $\Delta^{k-2}$ is coframed by $[df_{k-2} : du_{k-2} : dv_{k-2}]$, but since $p_{k-1} \in \lambda T_2$, it must have the form $p_{k-1}=(p_{k-2},l_{k-2})$ with $l_{k-2} \subseteq T_2(p_{k-2})$.  We must therefore have $T_2(p_{k-2})$ coframed by
\[[df_{k-2} : du_{k-2} : dv_{k-2}]= \left[1: \frac{du_{k-2}}{df_{k-2}} : \frac{dv_{k-2}}{df_{k-2}}\right]=[1 : u_{k-1} : v_{k-1}]\]
where $v_{k-1}=0$ and $u_{k-1}$ is not identically zero. Moreover, we see that $df_{k-1}=df_{k-2}$.  

Since $p_{k} \in \lambda T_2T_2$, the same argument shows that $T_2(p_{k-1})$ is coframed by
\[[df_{k-2} : du_{k-1} : dv_{k-1}]= \left[1: \frac{du_{k-1}}{df_{k-2}} : \frac{dv_{k-1}}{df_{k-2}}\right]=[1 : u_{k} : v_{k}]\]
where $v_{k}=0$ and $u_{k}$ is not identically zero. Moreover, we see that $df_k=df_{k-1}=df_{k-2}$.  

Now as an ansatz, suppose $T_2(p_k)$ indeed appears in $\Delta^k$.  Then it would have the form $ [df_{k-2} : du_{k} : dv_{k}]=[a : b : 0]$ with $a\neq0$.  Its projection one level down would have the form $ [df_{k-2} : du_{k-1} : dv_{k-1}]=[a : b : 0]$ with $a\neq0$.  We recognize this as $T_2(p_{k-1})$, which we know exists in $\Delta^{k-1}$.  Therefore $T_2(p_k)$ indeed exists as it is the prolongation of $T_2(p_{k-1})$, and our ansatz is justified.  

Finally, assume for sake of contradiction that  $T_1(p_k)$ appears in $\Delta^k$.  It would have the form $ [df_{k-2} : du_{k} : dv_{k}]=[a : 0 : b]$ with $a\neq0$. Its projection one level down would have the form $[df_{k-2} : du_{k-1} : dv_{k-1}]=[a : 0 : b]$ with $a\neq0$.  This forces $du_{k-1}=0$.  But we saw above that a local fiber coordinate at $p_{k-1}$ is $u_k=\frac{du_{k-1}}{df_{k-2}}$, and $u_k$ is not identically zero.  This contradiction disproves the existence of $T_1(p_k)$ in $\Delta^k$.

\subsubsection*{Case 2: $\omega=\lambda L_1T_2$}
Suppose $p_k\in M^k$ has RVT code $\omega=\lambda L_1T_2$.
From the discussion above, it suffices to prove that $T_2$ appears in $\Delta^k$, while $T_1$ does not.  Now the distribution $\Delta^{k}$ is coframed by $[df_{k} : du_{k} : dv_{k}]$.  Two levels down,  $\Delta^{k-2}$ is coframed by $[df_{k-2} : du_{k-2} : dv_{k-2}]$, but since $p_{k-1} \in \lambda L_1$, it must have the form $p_{k-1}=(p_{k-2},l_{k-2})$ with $l_{k-2} = L_1(p_{k-2})$.  We must therefore have $L_2(p_{k-2})$ coframed by
\[[df_{k-2} : du_{k-2} : dv_{k-2}]= \left[ \frac{df_{k-2}}{dv_{k-2}} : \frac{du_{k-2}}{dv_{k-2}}:1\right]=[ u_{k-1} : v_{k-1} : 1].\]
Moreover, we see that $df_{k-1}=dv_{k-2}$.  

Since $p_{k} \in \lambda T_2T_2$, the we can similarly see that $T_2(p_{k-1})$ is coframed by
\[[dv_{k-2} : du_{k-1} : dv_{k-1}]= \left[1: \frac{du_{k-1}}{dv_{k-2}} : \frac{dv_{k-1}}{dv_{k-2}}\right]=[1 : u_{k} : v_{k}]\]
where $v_{k}=0$ and $u_{k}$ is not identically zero. Moreover, we see that $df_k=dv_{k-2}$.  

Now as an ansatz, suppose $T_2(p_k)$ indeed appears in $\Delta^k$.  Then it would have the form $ [dv_{k-2} : du_{k} : dv_{k}]=[a : b : 0]$ with $a\neq0$.  Its projection one level down would have the form $ [dv_{k-2} : du_{k-1} : dv_{k-1}]=[a : b : 0]$ with $a\neq0$.  We recognize this as $T_2(p_{k-1})$, which we know exists in $\Delta^{k-1}$.  Therefore $T_2(p_k)$ indeed exists as it is the prolongation of $T_2(p_{k-1})$, and our ansatz is justified.  

Finally, assume for sake of contradiction that  $T_1(p_k)$ appears in $\Delta^k$.  It would have the form $ [dv_{k-2} : du_{k} : dv_{k}]=[a : 0 : b]$ with $a\neq0$. Its projection one level down would have the form $[dv_{k-2} : du_{k-1} : dv_{k-1}]=[a : 0 : b]$ with $a\neq0$.  This forces $du_{k-1}=0$.  But we saw above that a local fiber coordinate at $p_{k-1}$ is $u_k=\frac{du_{k-1}}{dv_{k-2}}$, and $u_k$ is not identically zero.  This contradiction disproves the existence of $T_1(p_k)$ in $\Delta^k$.

\subsection*{Acknowledgements}

The authors thank Richard Montgomery (Santa Cruz) for many useful conversations and remarks, and for introducing the authors to this subject.  Warm thanks also to Gary Kennedy (Ohio State) and Susan Colley (Oberlin) for continuing motivation.  We are also grateful to the referee for helpful comments and suggestions.

\bibliographystyle{amsplain}

\begin{thebibliography}{20}


\bibitem{arnold}
      \newblock V. I. Arnol'd,
      \newblock Simple singularities of curves,
      \newblock \emph{Proc. Steklov Inst. Math.}, \textbf{226} (1999), 20--28.

\bibitem{C} E. Cartan,
Sur l'\'{e}quivalence absolue de certains syst\`{e}mes d'\'{e}quations diff\'{e}rentielles et sur certaines familles de courbes, \emph{Bull. Soc. Math. France}, \textbf{XLII} (1914), 12--48.


\bibitem{castrothesis}
      \newblock A. Castro, 
      \newblock Chains and Monsters: from Cauchy-Riemann geometry to Semple towers and singular space curves,
      \newblock PhD thesis, 2010.

\bibitem{kennedy}
      \newblock A. Castro, S. Colley, G. Kennedy, and C. Shanbrom,
      \newblock A coarse stratification of the Monster tower,
      \newblock  arXiv:1606.07931 [math.AG].

\bibitem{castro3}
      \newblock A. Castro and W. Howard,
      \newblock A Monster tower approach to Goursat multi-flags,
      \newblock \emph{Differential Geom. Appl.}, \textbf{30} (2012), 405--427.
      
\bibitem{bridges}
      \newblock A. Castro, W. Howard, and C. Shanbrom,
      \newblock Bridges between subRiemannian geometry and algebraic geometry,
      \newblock \emph{ Proceedings of 10th AIMS Conference on Dynamical Systems, Differential Equations, and Applications},  (2015).
       
\bibitem{castro2}
      \newblock A. Castro, R. Montgomery, and appendix by W. Howard,
      \newblock Spatial curve singularities and the Monster/Semple tower,
      \newblock \emph{Israel J. Math.}, \textbf{192} (2012), 381--427.



\bibitem{GKR} A. Giaro, A. Kumpera, C. Ruiz,
Sur la lecture correcte d'un r\'{e}sultat d'\'{E}lie Cartan, \emph{C. R. Acad. Sci. Paris}, \textbf{287} (1978), 241--244.


\bibitem{J} F. Jean,
The car with N trailers: characterisation of the singular configurations,
\emph{ESAIM: Control, Optimisation, and Calculus of Variations}, \textbf{1} (1996), 241--266.



      
\bibitem{kumpera}
      \newblock A. Kumpera and J.L. Rubin,
      \newblock Multi-flag systems and ordinary differential equations,
      \newblock Nagoya Math. J., \textbf{166} (2002), 1--27. 
      

\bibitem{kushner}
A. Kushner, V. Lychagin, and V. Ruvtsov,
\emph{Contact geometry and nonlinear differential equations}
Cambridge University Press, Cambridge, UK, (2006).


\bibitem{respondek}  S.J. Li and W. Respondek
	\newblock The geometry, controllability, and flatness property of the $n$-bar system, 
	\newblock \emph{Internat. J. Control},
	\newblock\textbf{84}  (2011), 834--850.


\bibitem{LR} F. Luca and J.J. Risler,
The maximum of the degree of nonholonomy for the car with N trailers,
\emph{Proceedings of the 4th IFAC Symposium on Robot Control}, Capri, (1994), 165--170.

\bibitem{MZ2} R. Montgomery and M. Zhitomirskii,
Geometric Approach to Goursat Flags,
\emph{Ann. Inst. H. Poincar\'{e} - AN}, \textbf{18} (2001), 459--493.

\bibitem{MZ1} R. Montgomery and M. Zhitomirskii,
Points and Curves in the Monster Tower,
\emph{Memoirs of the AMS}, \textbf{956} (2010).


\bibitem{Mo1} P. Mormul,
Geometric classes of Goursat flags and their encoding by small growth vectors,
\emph{Central European J. Math.}, \textbf{2} (2004), 859--883.


\bibitem{Mo3} P. Mormul,
Multi-dimensional Cartan prolongation and special $k$-flags, \emph{Geometric Singularity Theory}, Banach Center Publications, \textbf{65} (2004), 157--178.

\bibitem{Mo2} P. Mormul,
Small growth vectors of the compactifications of the contact systems on $J^r(1,1)$, \emph{Contemporary Mathematics}, \textbf{569} (2012), 123--141.


\bibitem{morpel}
      \newblock P. Mormul and F. Pelletier,
      \newblock Special 2-flags in lengths not exceeding four: a study in strong nilpotency of distributions
      \newblock arXiv:1011.1763 [math.DG].

\bibitem{pelletier1}
      \newblock F. Pelletier and M. Slayman,
      \newblock Articulated arm and special multi-flags,
      \newblock \emph{J. Math. Sci. Adv. Appl.}, \textbf{8} (2011), 9--41.
      
  \bibitem{pelletier2}
      \newblock F. Pelletier and M. Slayman,
      \newblock Configurations of an articulated arm and singularities of special multi-flags,
      \newblock \emph{SIGMA}, \textbf{10} (2014).

\bibitem{respondek2}  W. Respondek,
	\newblock Symmetries and minimal flat outputs of nonlinear control systems, 
	\newblock \emph{New Trends in Nonlinear Dynamics and Control and their Applications},
	\newblock  Lecture Notes in Control and Information Science, \textbf{295}  (2004), 64--86.

  \bibitem{semple} J. Semple,
	\newblock Singularities of Space Algebraic Curves,
	\newblock \emph{Proceedings of the London Mathematical Society},
	\newblock\textbf{44}  (1938), 149--174.    

	      
\bibitem{Shanbrom} C. Shanbrom, 
	\newblock The Puiseux characteristic of a Goursat germ,
	\newblock \emph{J. Dynamical and Control Systems},
	\newblock\textbf{20}  (2014), 33--46.

   \end{thebibliography}

\end{document}